\documentclass{amsart}
\usepackage{amssymb,amsmath,latexsym,amsfonts,amsthm,mathrsfs,lineno}
\usepackage{enumerate}

\newtheorem{thm}{Theorem}
\newtheorem*{thm*}{Theorem}

\newtheorem*{claim}{Claim}
\newtheorem{claimm}{Claim}

\newcommand{\N}{\mathbb{N}}
\newcommand{\R}{\mathbb{R}}

\newcommand{\funct}[2]{#1 \longrightarrow #2}

\author{Lionel Nguyen Van Th\'e}

\address{Aix Marseille Univ, CNRS, Centrale Marseille, I2M UMR 7373, 13453 Marseille, France}

\email{lionel.nguyen-van-the@univ-amu.fr}

\title{On a problem of Specker about Euclidean representations of finite graphs}

\subjclass[2000]{Primary: 05C62}
\keywords{Finite graph, Euclidean representation, dimension}

\date{August, 2018}

\begin{document}

\maketitle

\dedicatory{\begin{center}This paper is dedicated to the memory of Ernst Specker.\end{center}}

\begin{abstract}

Say that a graph $G$ is \emph{representable in $\R ^n$} if there is a map $f$ from its vertex set into the Euclidean space $\R ^n$ such that $\| f(x) - f(x')\| = \| f(y) - f(y')\|$ iff $\{ x,x'\}$ and $\{ y, y'\}$ are both edges or both non-edges in $G$. The purpose of this note is to present the proof of the following result, due to Einhorn and Schoenberg in \cite{ES}: if $G$ finite is neither complete nor independent, then it is representable in $\R ^{|G|-2}$. A similar result also holds in the case of finite complete edge-colored graphs.

\end{abstract}

\section{Introduction}

Given a (simple and loopless) graph $G$ and a natural number $n \in \N$, say that $G$ is \emph{representable in $\R ^n$} if there is a map $f$ from the vertex set of $G$ (which we will also denote by $G$ in the sequel) into the Euclidean space $\R ^n$ such that $\| f(x) - f(x')\| = \| f(y) - f(y')\|$ iff $\{ x,x'\}$ and $\{ y, y'\}$ are both edges or both non-edges in $G$. Classical results about 2-distance sets in Euclidean spaces \cite{B} show that if $G$ is representable in $\R^n$, then $|G|\leq \binom{n+2}{2}$ where $|G|$ denotes of vertices of $G$. Equivalently: \[ \sqrt{|G|-\frac{3}{2}}-\frac{3}{4}\leq n \enspace . \]

On the other hand, it has been known for a long time\footnote{To my knowledge, this result appeared first in \cite{R} together with several other results about Euclidean representations of graphs. It is also a consequence of Schoenberg's theorem quoted below.} that every finite graph is representable in $\R ^{|G|-1}$. It is also clear that if $G$ is complete (i.e. all pairs carry an edge) or independent (i.e. no pair carries an edge), then $G$ is not representable in $\R ^{|G|-2}$ and dimension $|G|-1$ is necessary. But what about the converse? If $G$ is neither complete nor independent, is it representable in $\R ^{|G|-2}$? According to Maurice Pouzet, who mentions it in \cite{P} in connection to the famous Ulam reconstruction problem, this question was asked by Ernst Specker around 1972. Indeed, it was mentioned in 1973 during the conference in honor of the sixtieth birthday of Paul Erd\H os. Nobody could figure out the answer. However, it already existed at that time, and even had been published by Einhorn and Schoenberg in \cite{ES}.   

\begin{thm}
\label{thm:rep}
Let $G$ be a finite graph. Assume that $G$ is neither complete nor independent. Then $G$ is representable in $\R ^{|G|-2}$.
\end{thm}

The purpose of this note is to present the corresponding proof, which is only an elementary result in \cite{ES}. More generally, given a complete edge-colored graph $(G, \lambda)$ (a complete graph $G$ together with a map $\lambda : \funct{G^2}{\R}$ such that $\lambda(x,x)=0$ and $\lambda(y,x)=\lambda(x,y)$) and $n \in \N$, say that $G$ is \emph{representable in $\R ^n$} when there is a map $f : \funct{G}{\R^n}$ such that \[ \| f(x) - f(x')\| = \| f(y) - f(y')\| \ \ \textrm{iff} \ \ \lambda (x,x') = \lambda(y, y') \enspace .\] 

Again, known results about $k$-distance sets in Euclidean spaces \cite{BBS} show that if $(G, \lambda)$ is representable in $\R^n$, then $|G|\leq \binom{k+n}{k}$. On the other hand, every finite $(G, \lambda)$ is representable in $\R^{|G|-1}$ and if $\lambda$ takes only one value, then $(G, \lambda)$ is representable in $\R^{|G|-1}$ but not $\R^{|G|-2}$. But if $\lambda$ takes at least two values, representability in $\R^{|G|-2}$ is always guaranteed: 

\begin{thm}
\label{thm:repcol}
Let $(G, \lambda)$ be a complete edge-colored graph. Assume that $\lambda$ takes at least two different values. Then $G$ is representable in $\R ^{|G|-2}$. 
\end{thm}

Note that Theorem \ref{thm:rep} is a simple consequence of Theorem \ref{thm:repcol} when $\lambda$ takes at most two values. Theorem \ref{thm:repcol} is proved using the following well-known result due to Schoenberg, which provides a characterization of those complete edge-labelled graphs that appear as metric subspaces of some Euclidean space: For a complete edge-colored graph $(G, \lambda)$, with $G=\{ x_i:1\leq i \leq |G|\}$, and $\lambda$ with positive values, say that $(G, \lambda)$ is \emph{isometric} to a subset of $\R^{|G|-1}$ when there is a map $\varphi : G \rightarrow \R^{|G|-1}$ such that for every $i,j\leq |G|$, $\|\varphi(x_{i})-\varphi(x_{j})\|=\lambda(x_{i},x_{j})$. For a matrix $M=(m_{ij})_{1\leq i, j \leq n}$, define \[ Q_M = \max \left\{ \sum_{1\leq i,j \leq n} m_{ij}v_i v_j : \sum_{k=1} ^n v_k ^2 = 1 \ \ \textrm{and} \ \ \sum_{k=1} ^n v_k = 0 \right\} \enspace .\] 

\begin{thm*}[Schoenberg \cite{S}]
Let $(G, \lambda)$ be a complete edge-colored graph where $G=\{ x_i:1\leq i \leq |G|\}$ and $\lambda$ takes positive values. Let $M = (\lambda(x_i, x_j)^2)_{1\leq i,j \leq |G|}$. Then $(G, \lambda)$ is isometric to a subset of $\R^{|G|-1}$ iff $Q_M \leq 0$. In that case, the dimension of the affine space spanned by $(G, \lambda)$ is $(|G|-1)$ iff $Q_M < 0$.
\end{thm*}

A word of caution here: even though we are interested in $(G, \lambda)$, the matrix $M$ to be considered is the matrix of the \emph{squares} of the values of $\lambda$!   

\

The paper is organized as follows: For the sake of completeness, we start in Section \ref{section:Schoenberg} with a proof of Schoenberg's theorem. We continue in Section \ref{pf thm rep} with a proof of Theorem \ref{thm:rep}. The scheme of the proof is then reproduced in Section \ref{pf thm repcol} to prove Theorem \ref{thm:repcol}. Finally, Section \ref{story} tells the story behind this project.

\section{Proof of Schoenberg's theorem}

\label{section:Schoenberg}

We start with the first part of the theorem. Assume that $(G, \lambda)$ is isometric to a subset of $\R^{|G|-1}$, as witnessed by $\varphi : G \rightarrow \R^{|G|-1}$. Write $y_{i}=\varphi(x_{i})$. For $v_{1},..., v_{|G|-1}\in \R$, let $v=\sum_{i=1}^{|G|-1} v_{i}(y_{i}-y_{|G|})$. Then 
\begin{align*}
\|v\|^{2} &= \sum_{1\leq i,j\leq |G|-1}\langle y_{i}-y_{|G|},y_{j}-y_{|G|}\rangle v_{i}v_{j}\\
&= \sum_{1\leq i,j\leq |G|-1}\frac{1}{2}\left(\| y_{i}-y_{|G|}\|^{2}+\| y_{j}-y_{|G|}\|^{2}-\| y_{i}-y_{j}\|^{2}\right)v_{i}v_{j} 
\end{align*}

Summing over the three terms separately, we may write this as: $$ \|v\|^{2} = \sum_{i=1}^{|G|-1}v_{i}\cdot\sum_{i=1}^{|G|-1}\|y_{i}-y_{|G|}\|^{2}v_{i}-\frac{1}{2}\left(\sum_{1\leq i,j\leq |G|-1}\| y_{i}-y_{j}\|^{2}v_{i}v_{j}\right)$$

Setting $\displaystyle v_{|G|}=-\sum_{i=1}^{|G|-1}v_{i}$, we obtain $$ \|v\|^{2} = -\frac{1}{2}\|y_{|G|}-y_{|G|}\|^{2}v_{|G|}v_{|G|}-\sum_{i=1}^{|G|-1}\|y_{i}-y_{|G|}\|^{2}v_{i}v_{|G|}-\frac{1}{2}\left(\sum_{1\leq i,j\leq |G|-1}\| y_{i}-y_{j}\|^{2}v_{i}v_{j}\right)$$

i.e. $$ \|v\|^{2} = -\frac{1}{2}\sum_{1\leq i,j\leq |G|}\| y_{i}-y_{j}\|^{2}v_{i}v_{j}= -\frac{1}{2}\sum_{1\leq i,j\leq |G|}\lambda(x_{i},x_{j})^{2}v_{i}v_{j} \quad (*)$$

From this, it follows directly that $Q_{M}\leq 0$.  

Conversely, assume that $Q_{M}\leq 0$. Tracking back the previous computation, it follows that the quadratic form associated to the matrix $P$ is positive, where $$P=\left(\frac{1}{2}\left(\lambda(x_{i},x_{|G|})^{2}+\lambda(x_{j},x_{|G|})^{2}-\lambda(x_{i},x_{j})^{2}\right)\right)_{1\leq i,j\leq |G|-1}$$

Therefore, there exist vectors $z_{1},...,z_{|G|-1}\in \R^{|G|-1}$ such that for all $i,j\leq |G|-1$, $$ \langle z_{i},z_{j}\rangle = \frac{1}{2}\left(\lambda(x_{i},x_{|G|})^{2}+\lambda(x_{j},x_{|G|})^{2}-\lambda(x_{i},x_{j})^{2}\right)$$

Define now $y_{i}=z_{i}$ for $i\leq |G|-1$, and $y_{|G|}=0_{\R^{|G|-1}}$. We claim that the set $\{ z_{1},...,z_{|G|}\}$ is isometric to $(G, \lambda)$. Indeed, for $i, j\leq |G|-1$, $$ \| y_{i}-y_{|G|}\|^{2}=\|z_{i}\|^{2}=\langle z_{i},z_{j}\rangle=\lambda(x_{i},x_{|G|})^{2}$$ 

and $$\| y_{i}-y_{j}\|^{2}=\|z_{i}-z_{j}\|^{2}=\|z_{i}\|^{2}+\| z_{j}\|^{2}-2\langle z_{i},z_{j}\rangle$$

Now, replacing $\langle z_{i},z_{j}\rangle$ by $\frac{1}{2}\left(\lambda(x_{i},x_{|G|})^{2}+\lambda(x_{j},x_{|G|})^{2}-\lambda(x_{i},x_{j})^{2}\right)$ leads to 
$$\| y_{i}-y_{j}\|^{2}=\lambda(x_{i},x_{j})^{2}$$

This finishes the proof of the first part of Schoenberg's theorem. For the second part, consider the vector $v$ defined as previously. The dimension of the affine space spanned by $(G, \lambda)$ is $(|G|-1)$ iff $\|v\|^{2}>0$ whenever not all the $v_{i}$'s are zero. By $(*)$ this means exactly that $Q_{M}<0$.

\section{Proof of Theorem \ref{thm:rep}}

\label{pf thm rep}

The geometric idea behind the proof is elementary: Start with a equilateral metric space on $|G|$ points in $\R^{|G|-1}$, where all distances are equal to $1$. Perturbing certain lengths to some number $\alpha\approx 1$, and other lengths to some other number $\beta\approx 1$, this can be arranged to become a Euclidean representation of $G$. When $G$ complete or independent, this is of course always the case regardless of the choice for $\alpha$ and $\beta$, but when this is not so, some other choice $(\alpha_{0}, \beta_{0})$ makes this configuration non-metric. Therefore, when varying continuously from $(1,1)$ to $(\alpha_{0}, \beta_{0})$, the configuration gets continuously deformed until reaching some point where it stops being Euclidean. This geometric obstruction materializes by a non-trivial affine relationship between the points of the configuration, and hence a drop in the dimension of the embedding of the representation, which becomes at most $|G|-2$. 

Let us now proceed with the detailed proof. Let $G$ be a finite graph that is neither complete nor independent. Enumerate the vertices of $G=\{ x_k:1\leq k \leq |G|\}$ and let $M_G=(m_{ij})_{1\leq i,j \leq |G|}$ denote the adjacency matrix of $G$ with respect to this enumeration, ie: 

\begin{displaymath}
m_{ij} = \left \{ \begin{array}{l}
 1 \ \ \textrm{if $\{x_i, x_j\}$ is an edge in $G$}, \\
 0 \ \ \textrm{otherwise}.
 \end{array} \right.
\end{displaymath}

Let $\overline{M}$ be the adjacency matrix of the complement of $G$ (the graph obtained from $G$ by changing all the edges between different vertices into non-edges and vice-versa). For $\alpha, \beta > 0$, let \[ M(\alpha, \beta) = \alpha M + \beta \overline{M} \enspace . \]

Denoting $M(\alpha, \beta) = (m^{\alpha \beta}_{ij})_{1\leq i,j \leq |G|}$, say that $M(\alpha, \beta)$ \emph{codes a representation of $G$ in $\R ^{|G|-1}$} when the complete edge-colored graph $(G, d)$, with $d (x_i, x_j)= m^{\alpha \beta}_{ij}$, is isometric to a subset of $\R ^{|G|-1}$. According to Schoenberg's theorem, we need to show that there are $\alpha \neq \beta > 0$ such that $Q_{M(\alpha^2, \beta^2)}=0$. (Recall that  Schoenberg's theorem relates to the matrix of the squares of $d$, hence the appearance of $\alpha^{2}$ and $\beta^{2}$ in the preceding expression.)

\begin{claimm}

\label{claim1}

There are $\alpha_0, \beta_0 > 0$ such that $Q_{M(\alpha_0 ^2, \beta_0 ^2)} > 0$.  

\end{claimm}

\begin{proof}

Assume towards a contradiction that $Q_{M(\alpha ^2, \beta ^2)} \leq 0$ for all $\alpha, \beta > 0$. We show that $G$ is complete or independent. Indeed, first take $\alpha, \beta > 0$ such that $2\alpha < \beta$. Since $Q_{M(\alpha ^2, \beta ^2)} \leq 0$, Schoenberg's theorem guarantees that $M(\alpha, \beta)$ codes a representation of $G$ in $\mathbb{R}^{|G|-1}$ and by triangle inequality, no triangle with two sides of length $\alpha$ and one side of length $\beta$ appears in this representation. Therefore, $G$ does not contain the graph $H$ drawn in Figure \ref{fig:H}.

\begin{figure}[h]
\setlength{\unitlength}{0.5mm}
\begin{picture}(20,30)(0,0)

\put(10,25){\circle*{2}}
\put(0,5){\circle*{2}}
\put(20,5){\circle*{2}}
\put(9,24){\line(-1,-2){9}}
\put(11,24){\line(1,-2){9}}

\end{picture}

\caption{The graph $H$.}

\label{fig:H}

\end{figure}
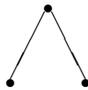

Similarly, choosing $2\beta < \alpha$, no triangle with one side of length $\alpha$ and two sides of length $\beta$ appears in the representation coded by $M(\alpha, \beta)$. Therefore, $G$ does not contain the graph $K$ depicted in Figure \ref{fig:K}. 

\begin{figure}[h]
\setlength{\unitlength}{0.5mm}
\begin{picture}(20,30)(0,0)

\put(10,25){\circle*{2}}
\put(0,5){\circle*{2}}
\put(20,5){\circle*{2}}
\put(1,5){\line(1,0){18}}

\end{picture}

\caption{The graph $K$.}
\label{fig:K}

\end{figure}
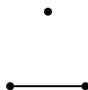

It follows that $G$ is complete or independent, a contradiction. \end{proof}

\begin{claimm}
\label{claim2} 
The map $M \mapsto Q_M$ is continuous ($n \times n$ matrices are seen as elements of $\R^{n^2}$ equipped with the standard topology). 
\end{claimm}

\begin{proof}
Since the topology of $\R^{n^2}$ is the topology induced by the $\ell_1$ norm (ie $\| M \| = \sum_{1\leq i, j \leq n} |m_{ij}|$), it is enough to show that $|Q_M - Q_N|\leq \|M - N\|$. This is done by observing that whenever $\sum_{k=1} ^n v_k ^2 = 1$, we have  
\begin{align*} 
\left| \sum_{1\leq i,j \leq n} m_{ij}v_i v_j - \sum_{1\leq i,j \leq n} n_{ij}v_i v_j \right| & \leq \sum_{1\leq i,j \leq n} |m_{ij} - n_{ij}||v_i v_j| \\ 
& \leq \sum_{1\leq i,j \leq n} |m_{ij} - n_{ij}| \\
& \leq \| M-N\| \enspace .
\end{align*}

Therefore 

\begin{align*}
Q_M & = \max \left\{ \sum_{1\leq i,j \leq n} m_{ij}v_i v_j : \sum_{k=1} ^n v_k ^2 = 1 \ \ \textrm{and} \ \ \sum_{k=1} ^n v_k = 0 \right\} \\ 
& \leq \max \left\{ \sum_{1\leq i,j \leq n} n_{ij}v_i v_j + \| M-N\| : \sum_{k=1} ^n v_k ^2 = 1 \ \ \textrm{and} \ \ \sum_{k=1} ^n v_k = 0 \right\} \\
& \leq Q_N + \| M - N\| \enspace .
\end{align*}

Hence, $Q_M - Q_N \leq \| M - N\|$ and by symmetry, $Q_N - Q_M \leq \| M - N\|$. It follows that $|Q_M - Q_N|\leq \|M - N\|$. \end{proof}

By Claim \ref{claim1} pick $\alpha_0, \beta_0 > 0$ such that $Q_{M(\alpha_0 ^2, \beta_0 ^2)} > 0$. Note that without loss of generality, we may assume that $\alpha_0 \neq \beta_0$. This is because continuity of the map $M \mapsto Q_M$ proved in Claim \ref{claim2} implies continuity of $(\alpha, \beta) \mapsto Q_{M(\alpha, \beta)}$. For $t \in [0,1]$, consider the matrix \[ M\left(1+t(\alpha_0 ^2 - 1), 1+t(\beta_0 ^2 - 1) \right) \enspace .\]

It defines a continuous curve from $M(1,1)$ to $M(\alpha_0 ^2 , \beta_0 ^2)$, and the map \[ \psi : t \mapsto Q_{M\left(1+t(\alpha_0 ^2 - 1), 1+t(\beta_0 ^2 - 1)\right)}\] is continuous on $[0,1]$. Observe that $M(1,1)$ codes the equilateral metric space on $|G|$ points where all the distances are equal to one. This metric space is Euclidean and spans an affine space of dimention $|G|-1$, therefore $\psi(0) = Q_{M(1,1)} < 0 $. Observe on the other hand that $\psi(1) = Q _{M(\alpha_0 ^2, \beta_0 ^2)}>0$. So by the intermediate value theorem, there is $\tau \in (0,1)$ such that $\psi(\tau)=0$. That means \[Q_{M\left(1+\tau (\alpha_0 ^2 - 1), 1+\tau(\beta_0 ^2 - 1) \right)}=0 \enspace .\] 

So set $\alpha = \sqrt{1+\tau(\alpha_0 ^2 - 1)}$ and $\beta = \sqrt{1+\tau(\beta_0 ^2 - 1)}$. Then $\alpha \neq \beta > 0$ and $M\left(\alpha, \beta\right)$ codes a representation of $G$ in $\R^{|G|-2}$. \qed

\

As the very last step of the preceding proof is non-constructive, it is natural to ask whether an exact computation of $(\alpha, \beta)$ coding a representation of $G$ in $\R^{|G|-2}$ could be performed, and whether this could result in a representation of lower dimension. And indeed, this is the case, as shown by Roy in \cite{Ro}. The results of this paper show that the value $\alpha = \tau/(\tau+1)$, $\beta=1-\alpha$ usually works, where $\tau$ is the smallest eigenvalue of $M$, and exhibit the exact computation of the minimal dimension into which $G$ can be represented.

\section{Proof of Theorem \ref{thm:repcol}}

\label{pf thm repcol}

The proof follows exactly the same pattern as the proof of Theorem \ref{thm:rep} so we only emphasize the ideas. Let $(G, \lambda)$ be a complete colored graph where $\lambda$ has range $\{l_1,\ldots , l_p\}$ of size at least two. Enumerate the vertices of $G=\{ v_k:1\leq k \leq |G|\}$ and let $M_i$ denote the adjacency matrix of the graph obtained from $G$ by keeping only the edges with color $l_i$. For $\alpha_1, \ldots, \alpha_p > 0$, let \[ M(\alpha_1, \ldots, \alpha_p) = \sum _{i=1} ^p \alpha_i M_i \enspace . \]

According to Schoenberg's theorem, we need to show that there are distinct $\alpha_1, \ldots, \alpha _p > 0$ such that $Q_{M(\alpha_1 ^2, \ldots , \alpha_p ^2)}=0$. 

\begin{claim}

There are $a_1, \ldots, a_p > 0$ such that $Q_{M(a_1 ^2, \ldots , a_p ^2)} > 0$. 

\end{claim}

\begin{proof}

Suppose not. Then $Q_{M(\alpha_1 ^2, \ldots , \alpha_p ^2)} \leq 0$ for all $\alpha_1, \ldots, \alpha _p > 0$. Varying the coefficients $\alpha_1, \ldots, \alpha _p$ and taking, turn by turn, $\alpha_i$ much larger than all the other coefficients, triangle inequality in the corresponding representations shows that all the triangles in $(G, \lambda)$ must have all their egdes of the same color. Therefore, $\lambda$ only takes one value, a contradiction. \end{proof}

So pick $a_1, \ldots, a_p > 0$ such that $Q_{M(a_1 ^2, \ldots , a_p ^2)} > 0$. Note that the continuity of the map $M \mapsto Q_M$ (Claim \ref{claim2}) guarantees that without loss of generality, we may assume that all the $a_i$'s are distinct. For $t \in [0,1]$, consider the matrix \[M\left(1+t(a_1 ^2 - 1), \ldots, 1+t(a_p ^2 - 1) \right) \enspace .\] 

It defines a continuous curve from $M(1,\ldots, 1)$ to $M(a_1 ^2 , \ldots, a_p ^2)$, and the map \[ \psi : t \mapsto Q_{M\left(1+t(a_1 ^2 - 1), \ldots, 1+t(a_p ^2 - 1)\right)}\] is continuous on $[0,1]$. Since $M(1,\ldots, 1)$ codes a Euclidean metric space that spans an affine space of dimention $|G|-1$, we have $\psi(0) = Q_{M(1,\ldots, 1)} < 0 $. On the other hand, $\psi(1) = Q _{M(a_1 ^2, \ldots, a_p ^2)}>0$. So by the intermediate value theorem, there is $\tau \in (0,1)$ such that $\psi(\tau)=0$. That means \[Q_{M\left(1+\tau (a_1 ^2 - 1), \ldots, 1+\tau(a_p ^2 - 1) \right)}=0 \enspace .\] 

So for $1\leq i \leq p$, set $\alpha _i = \sqrt{1+\tau(a_i ^2 - 1)}$. Then all the $\alpha_i$'s are $> 0$ and distinct, and $M\left(\alpha_1, \ldots, \alpha_p\right)$ codes a representation of $G$ in $\R^{|G|-2}$. \qed

\section{Afterword}

\label{story}

I first heard about the problem of Euclidean representation of finite graphs in 2004, when meeting Maurice Pouzet, who in turn had heard it from Ernst Specker around 1972. Pouzet had been advertizing the problem since then, but had never heard any progress about it. The situation had not changed when he visited Claude Laflamme, Norbert Sauer, Robert Woodrow and I in Calgary, in 2008. He mentioned the problem again. This is when the solution of the present paper was found, written down, and submitted.  

Around the same time, I realized that another colleague, Ilijas Farah, was in touch with Ernst Specker, who was still living in Zurich. Therefore, when I had the opportunity to visit Switzerland, I contacted him. Very kindly, he invited me to meet him. Of course, I was very excited to ask him how he had been led to this remarkable problem. He laughed out loud when he heard the question: there was no motivation at all! He had just noticed that he could prove the result by hand for small graphs, but that he could not do it in general. This is how the problem made it to Hungary and to Paul Erd\H os' birthday conference in 1973 (which was taking place close to Lake Balaton, in Hungary). It could be that Pouzet heard about it there, but this is not completely clear (even to Pouzet himself). 

In October 2010, the journal where I had submitted my paper wrote back. I was surprised to read that a reference from 1966 had been provided by the referee, and completely solved the problem! In fact, the problem itself is not even mentioned explicitly, and its solution only appears as one of the elementary results at the beginning (Lemma 2). And sure enough, the proof I had was nothing different than what appeared there. So of course, my paper was not published. However, as mentioned previously, something else was: Aidan Roy, who worked at the University of Calgary at the same time as I did, managed to come up with the exact minimal dimension that is necessary to embed a given graph. The proof is much more sophisticated than those contained here and uses spectral graph theory. It can be found in \cite{Ro}. 

The last act of this little play takes place in April 2012. Ernst Specker had just passed away, on December 10th, 2011. Alain Valette, by whom I heard this sad news, suggested that my little unpublished paper would be in place in Expositiones Mathematicae. It took me a few more years to take the time to make the appropriate modifications, but I am sincerely glad and honored to see the curtain fall that way.  

\

\textbf{Acknowledgements}: This not so short story contains many characters, without whom it would not really have been a story. I am particularly indebted to Maurice Pouzet, without whom the problem would not have survived for so long; Alain Valette, without whom nothing would have been  published about it; and of course, Ernst Specker, without whom nothing would have happened.


\begin{thebibliography}{BBS}

\bibitem[BBS83]{BBS}
E. Bannai, E. Bannai and D. Stanton, An upper bound for the cardinality of an $s$-distance subset in real Euclidean space. II, \emph{Combinatorica}, 3 (2), 1983, 147--152. 

\bibitem[B81]{B}
A. Blokhuis, A new upper bound for the cardinality of 2-distance sets in Euclidean space, \emph{North-Holland Math. Stud.}, 87 [Special Issue: Convexity and graph theory, Jerusalem, Isra\"el, 1981], 1984, 65--66. 

\bibitem[ES66]{ES}
S.J. Einhorn and I.J. Schoenberg, On Euclidean sets having only two distances between points. I, II, 
\emph{Nederl. Akad. Wet., Proc., Ser. A} 69, 479-488, 1966, 489--504.


\bibitem[P79]{P}
M. Pouzet, Sur le probl\`eme de Ulam, \emph{J. Combin. Theory Ser. B}, 27 (3), 1979, 231--236, \emph{in French}.

\bibitem[R84]{R}
F. Reverdy, Repr\'esentation des graphes dans les espaces euclidiens et probl\`emes de repr\'esentation, M\'emoire de D.E.A. , Universit\'e Lyon 1, 1984 (French). 

\bibitem[Ro10]{Ro}
A. Roy, Minimal Euclidean representations of graphs, \emph{Discrete Mathematics}, 310 (4), 2010, 727--733. (arxiv.org/abs/0812.3707).

\bibitem[S38]{S}
I. J. Schoenberg, Metric spaces and positive definite functions, \emph{Trans. Amer. Math. Soc.}, 44 (3), 1938, 522--536.
\end{thebibliography}
\end{document}